\documentclass[12pt]{article}

\usepackage{amssymb,amsmath}
\usepackage{theorem}
\usepackage{graphicx}

\newcommand{\C}{{\mathbb C}}

\newcommand{\Z}{{\mathbb Z}}
\newcommand{\fsl}{{\mathfrak{sl}}}
\newcommand{\g}{{\mathfrak g}}

\newcommand{\cM}{{\cal M}}
\newcommand{\cPM}{{\cal PM}}
\newcommand{\cDM}{{\cal DM}}

\newcommand{\G}{{\Gamma}}
\newcommand{\tG}{\widetilde{\Gamma}}

\newtheorem{theorem}{Theorem}[section]

\newtheorem{conjecture}[theorem]{Conjecture}
\newtheorem{proposition}[theorem]{Proposition}
\newtheorem{corollary}[theorem]{Corollary}

{\theorembodyfont{\rmfamily}    
\newtheorem{remark}[theorem]{Remark}
\newtheorem{example}[theorem]{Example}
\newtheorem{definition}[theorem]{Definition}

}

\newcommand\mdskip{\ifnum\medskipamount>\lastskip
      \vskip-\lastskip\vskip\medskipamount\fi}
\newcommand\QED{\ifmmode\eqno\square\else
{\parfillskip0pt\hfil$\square$\par}\mdskip\fi}

\title{On a weight system conjecturally related to $\fsl_2$}
\author{E.~Kulakova, S.~Lando\thanks{Partly supported by the Russian Foundation for Basic Research Grant
no.~13-01-00383a}, T.~Mukhutdinova\thanks{Partly supported by the Russian Foundation for Basic Research Grant
no.~13-01-00383a}, G.~Rybnikov\\
National Research University Higher School of Economics}
\date{}

\begin{document}

\maketitle
\begin{abstract}
We introduce a new series~$R_k$, $k=2,3,4,\dots$, of integer valued weight systems.
The value of the weight system~$R_k$ on a chord diagram is
a signed number of cycles of even length~$2k$ in the intersection
graph of the diagram. We show that this value depends on the intersection
graph only. We check that for small orders of the diagrams,
the value of the weight system~$R_k$ on a diagram of order exactly~$2k$
coincides with the coefficient of~$c^k$ in the value of the $\fsl_2$-weight
system on the projection of the diagram to primitive elements.
\end{abstract}

\section{Introduction}

Below, we use standard notions from the theory of finite order knot invariants;
see, e.g.~\cite{CDM12} or~\cite{LZ04}.

\begin{figure}[ht]
\center{\includegraphics[width=.7\linewidth]{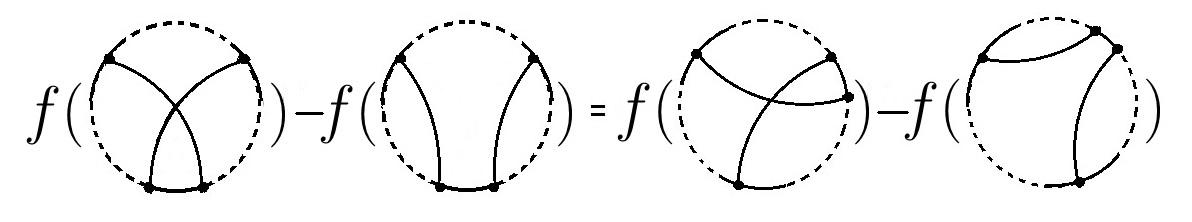}}
		\caption{$4$-term relation}
		\label{fourtermrelation}
\end{figure}

A {\it chord diagram\/} of order~$n$ is an oriented circle endowed with~$2n$
pairwise distinct points split into~$n$ disjoint pairs, considered up to
orientation-preserving diffeomorphisms of the circle.
A~{\it weight system\/} is a function on chord diagrams satisfying the $4$-term
relation; see Fig.~\ref{fourtermrelation}. For a chord diagram~$d$ with two chords~$A$
and~$B$ having neighboring ends, as in Fig.~\ref{fourtermrelation}, we will write this relation as
$f(d)-f(d'_{AB})=f(\tilde d_{AB})-f(\tilde d'_{AB})$.

In figures, the outer circle of the chord diagram is always assumed to be
oriented counterclockwise. Dashed arcs may contain ends of arbitrary sets
of chords, same for all the four terms in the picture.

To each weight system, a finite order knot invariant can be associated in a canonical way,
which makes studying weight systems an important part of knot theory.
There is a number of approaches to constructing weight systems. In particular,
weight systems can be constructed from semisimple Lie algebras,
although the result is complicated.
The present paper has been motivated by an aspiration for understanding the weight
system corresponding to the simplest nontrivial case of the Lie algebra~$\fsl_2$.

The {\it intersection graph\/} $\gamma(d)$ of a chord diagram~$d$ is the simple graph whose
vertices are in one-to-one correspondence with the chords in~$d$,
and two vertices are connected by an edge iff the corresponding chords intersect
one another. The {\it $4$-term relation for graphs}, introduced in~\cite{L00},
is the graph counterpart of the $4$-term relation for chord diagrams.
It is defined as follows.

Denote by $V(\G)$ the set of vertices of a graph $\G$ and by
$E(\G)$ the set of its edges.
Let us associate to each ordered pair of (distinct)
vertices $A,B\in V(\G)$ of a graph
$\G$ two other graphs $\G_{AB}'$ and $\tG_{AB}$.

The graph $\G_{AB}'$ is obtained from $\G$ by erasing the edge $AB\in E(\G)$
in the case that this edge exists, and by adding the edge otherwise.
In other words, we simply change the adjacency between the vertices
$A$ and $B$ in $\Gamma$. This operation is an analogue of edge deletion,
but we prefer to formulate it in a slightly more symmetric way.

The graph $\tG_{AB}$ is obtained from $\G$ in the following way. For any vertex
$C\in V(\G)\setminus\{A,B\}$ we switch its adjacency with $A$
to the opposite one if $C$ is joined with $B$,
and we do nothing otherwise. All other edges do not change.
Note that the graph $\tG_{AB}$ depends not only on the pair $(A,B)$,
but on the order of vertices  in the pair as well.

\begin{definition}[4-invariant]\rm
A graph invariant~$f$ is a
\index{4-invariant!of a graph}%
{\it $4$-invariant}\/ if it satisfies
the {\it $4$-term relation}
\begin{equation}\label{7eq:21}
f(\G)-f(\G_{AB}')=f(\tG_{AB})-f(\tG_{AB}')
\end{equation}
for each graph $\G$ and for any pair $A,B\in V(\G)$ of its vertices.
\end{definition}

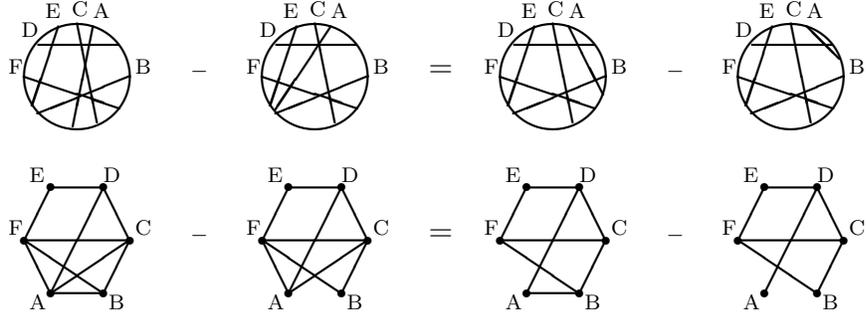
\begin{figure}[ht]
\begin{center}
\begin{picture}(200,100)(60,0)
\thicklines
\multiput(12,72)(90,0){4}{\circle{40}}
\multiput(55,72)(180,0){2}{{\scriptsize $-$}}
\put(145,72){$=$}
\multiput(5,91)(90,0){4}{\line(-1,-3){10}}
\multiput(0,94)(90,0){4}{{\scriptsize E}}
\multiput(-3,84)(90,0){4}{\line(1,0){31}}
\multiput(-9,87)(90,0){4}{{\scriptsize D}}
\multiput(-8,72)(90,0){4}{\line(3,-1){36}}
\multiput(-14,73)(90,0){4}{{\scriptsize F}}
\multiput(32,72)(90,0){4}{\line(-5,-2){35}}
\multiput(34,73)(90,0){4}{{\scriptsize B}}
\multiput(12,92)(90,0){4}{\line(1,-5){7.5}}
\multiput(10,95)(90,0){4}{{\scriptsize C}}
\put(18,91){\line(-1,-5){7.6}}
\multiput(18,94)(90,0){4}{{\scriptsize A}}
\put(108,91){\line(-2,-3){20.9}}
\put(198,91){\line(1,-2){13}}
\put(288,91){\line(1,-1){12.3}}

\multiput(2,30)(90,0){4}{\circle*{3}}
\multiput(22,30)(90,0){4}{\circle*{3}}
\multiput(2,-10)(90,0){4}{\circle*{3}}
\multiput(22,-10)(90,0){4}{\circle*{3}}
\multiput(-8,10)(90,0){4}{\circle*{3}}
\multiput(32,10)(90,0){4}{\circle*{3}}
\multiput(55,10)(180,0){2}{{\scriptsize $-$}}
\put(145,10){$=$}
\multiput(-6,-16)(90,0){4}{{\scriptsize A}}
\multiput(24,-16)(90,0){4}{{\scriptsize B}}
\multiput(-6,32)(90,0){4}{{\scriptsize E}}
\multiput(22,32)(90,0){4}{{\scriptsize D}}
\multiput(-14,12)(90,0){4}{{\scriptsize F}}
\multiput(34,12)(90,0){4}{{\scriptsize C}}

\multiput(-8,10)(90,0){4}{\line(1,0){39}}
\multiput(-8,10)(90,0){4}{\line(1,2){10}}
\multiput(-8,10)(90,0){2}{\line(1,-2){10}}
\multiput(-8,10)(90,0){4}{\line(3,-2){30}}
\multiput(32,10)(90,0){2}{\line(-3,-2){30}}
\multiput(32,10)(90,0){4}{\line(-1,2){10}}
\multiput(32,10)(90,0){4}{\line(-1,-2){10}}
\multiput(2,30)(90,0){4}{\line(1,0){20}}
\multiput(2,-10)(180,0){2}{\line(1,0){20}}
\multiput(2,-10)(90,0){4}{\line(1,2){20}}

\end{picture}
\end{center}
		\caption{A $4$-term relation for chord diagrams and the corresponding
intersection graphs}
		\label{fourtermrelationfor graphs}
\end{figure}

Figure~\ref{fourtermrelationfor graphs} shows the $4$-term relation for intersection graphs
corresponding to a sample $4$-term relation for chord diagrams.

Any graph invariant satisfying the $4$-term relation (that is, a $4$-invariant)
determines a weight system~\cite{L00}: the value of this weight system on a chord diagram
is set to be the value of the $4$-invariant on the intersection graph of the diagram.

Let~$d$ be a chord diagram. Denote by~$E_l(d)$
the number of circuits of length~$l$ ($=$~edge $l$-gons) in the intersection graph of~$d$.
It is well known (see~\cite{LZ04}, Exercise~6.4.10) that, for~$l\ge4$,
the parity of the number~$E_l(d)$ is a weight system with values in~$\Z/2\Z$.
However, no integer-valued weight system with the same parity has been known.

In this paper, we construct, for each $k=2,3,4,\dots$, an integer-valued weight system
$R_k$ whose value on each chord diagram has the same parity as $E_{2k}$, $R_k(d)\equiv E_{2k}(d)\mod~2$
for any chord diagram~$d$.
This weight system counts circuits of length~$2k$ in the intersection graph
of the diagram with signs depending on the mutual position of the corresponding chords.

For completeness, let us recall the proof of the fact mentioned above.
It is worth to be compared with the proof for the integer-valued analogue below.

\begin{proposition}\label{pmod2}
The value~$E_l\mod~2$ is a weight system for each $l\ge4$.
\end{proposition}

{\bf Proof.}
The value $E_l(d)$ obviously depends on the intersection graph
of a chord diagram rather than on the diagram itself. Let us extend
the function~$E_l$ to arbitrary graphs (that are not necessarily intersection graphs)
in the obvious way: let $E_l(\Gamma)$ be the number of edge $l$-gons in~$\Gamma$ having~$l$ pairwise
distinct vertices. We are going to prove that $E_l(\Gamma)\mod~2$
satisfies the $4$-term relation for graphs; the proposition then follows.

Suppose~$\Gamma$ contains an edge~$AB$. Then for the two terms on
the left-hand side of the $4$-term relation  we have that
$E_l(\Gamma)-E_l(\Gamma'_{AB})$
is the number of edge $l$-gons in~$\Gamma$ passing through the edge $AB$. Similarly,
for the right-hand side,
$E_l(\widetilde\Gamma)-E_l(\widetilde\Gamma'_{AB})$
is the number of edge $l$-gons in
$\widetilde\Gamma_{AB}$ passing through $AB$.
All the $l$-gons in~$\Gamma$
 passing through $AB$ contain a chain $CABD$ and split
into three disjoint classes according to the adjacency of the vertices~$C$ and~$D$
to~$A$ and~$B$:
\begin{enumerate}
\item the vertex~$C$ is adjacent to~$B$ and~$D$ is adjacent to~$A$;
\item the vertex~$C$ is adjacent to~$B$, but~$D$ is not adjacent to~$A$;
\item the vertex C is not adjacent to B.
\end{enumerate}
All edge $l$-gons in~$\widetilde\Gamma_{AB}$ passing through the four points $A,B,C,D$ admit a
similar classification.

\begin{example}
The three edge $4$-gons in the leftmost graph in Fig.~\ref{fourtermrelationfor graphs} passing through the edge~$AB$
split in the three classes in the following way:
\begin{itemize}
\item the edge quadrangles $FABC$ and $CABF$ belong to the first class;
\item there are no edge quadrangles belonging to the second class;
\item the edge quadrangle $DABC$ belongs to the third class.
\end{itemize}
\end{example}

Now, the $l$-gons in~$\Gamma$ belonging to the second class are in one-to-one
correspondence with the $l$-gons in~$\widetilde\Gamma_{AB}$ containing the path $CBAD$. The
edge $l$-gons of the third kind are the same in both graphs~$\Gamma$ and~$\widetilde\Gamma_{AB}$. And,
finally, the edge $l$-gons of the first kind in each of the two graphs come in
pairs: the chain $CABD$ can be replaced with the chain $CBAD$. Hence, the
number of edge $l$-gons of the first type is even for each of the two graphs, and
the required assertion follows. \QED

The authors are grateful to the participants of the seminar ``Combinatorics of finite order knot invariants''
at the Department of mathematics NRU HSE for useful discussions. The authors would also like to express
their gratitude to the referees for careful proofreading and valuable suggestions.

\section{Definition of the weight systems~$R_k$}\label{sd}

 \subsection{Positively and negatively oriented $2k$-gons}

\begin{figure}[ht]
\center{\includegraphics[width=.5\linewidth]{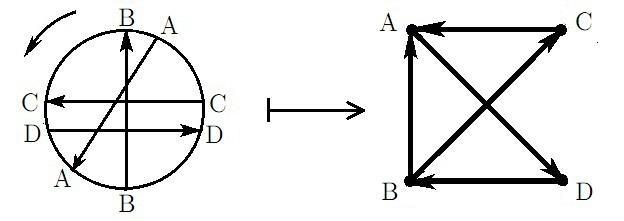}}
		\caption{An oriented chord diagram and the corresponding directed intersection graph;
the $4$-circuit $ADBC$ is oriented positively}
		\label{feo}
\end{figure}

 In order to define the weight system~$R_k$, let us take a chord diagram~$d$
 and choose an arbitrary orientation of its chords. This orientation induces
 an orientation of the edges of the intersection graph $\gamma(d)$ in the following
 way. We orient an edge~$AB$ from~$A$ to~$B$ if the beginning of the chord~$B$
 belongs to the arc of the outer circle of~$d$ which starts
 at the beginning of~$A$ and goes in the positive direction to the end of~$A$; see Fig.~\ref{feo}.
We denote the directed intersection graph of an oriented chord diagram~$d$ by~$\gamma(d)$,
like in the case of ordinary intersection graph of a non-oriented chord diagram, since
this convention causes no misunderstanding. Oriented edges in a directed graph will also be called {\it arrows}.

 We say that a circuit of even length $l=2k$ in a directed graph
 is {\it positively oriented\/} (or its sign is~$+$) if the number of arrows in this circuit
 oriented in either direction is even; otherwise the circuit is {\it negatively oriented}
 (or has the sign~$-$).
 Since the total number of arrows is even, the sign is independent of the choice of
 the direction of the circuit. For example, all the arrows in the $4$-circuit~$ADBC$
 in the directed intersection graph in Fig.~\ref{feo} are oriented in the same direction,
 meaning it is positive.

 \begin{definition}
 The function~$R_k$ takes a chord diagram~$d$ to the difference between the number
 of positively and negatively oriented $2k$-gons in the directed intersection graph
 $\gamma(d)$ for arbitrarily chosen orientation of the chords in~$d$.
 \end{definition}

Now we are going to show that this definition makes sense.

\begin{proposition}
The function~$R_k$ is well defined, that is,
its value on a chord diagram does not depend on the chosen orientation of the
chords.
\end{proposition}

{\bf Proof.} It suffices to check that if we switch the direction of a single chord
in a chord diagram~$d$ to the opposite one, then the sign of each circuit remains the same.
Changing the direction of a chord~$A$ in~$d$ leads to switching directions of all the arrows
in~$\gamma(d)$ incident to the vertex~$A$, preserving the directions of all other arrows.
Hence, the sign of any circuit not containing the vertex~$A$ remains the same. Any circuit
passing through~$A$ contains exactly two arrows incident to~$A$. After switching
the direction of both of them, the sign also remains the same. \QED

\begin{figure}[ht]
\begin{center}
\begin{picture}(200,40)(50,0)
\thicklines
\put(12,12){\circle{40}}
\put(45,10){$\mapsto$}
\put(-2,25){\vector(1,-1){28}}
\put(-2,-1){\vector(1,1){28}}
\put(12,32){\vector(0,-1){40}}
\put(-8,12){\vector(1,0){40}}

\put(9,34){{\scriptsize A}}
\put(-9,27){{\scriptsize B}}
\put(-13,-4){{\scriptsize D}}
\put(-16,10){{\scriptsize C}}

\put(70,-5){\circle*{3}}
\put(70,25){\circle*{3}}
\put(100,-5){\circle*{3}}
\put(100,25){\circle*{3}}
\put(70,-5){\vector(1,1){30}}
\put(70,-5){\vector(1,0){30}}
\put(70,-5){\vector(0,1){30}}
\put(100,-5){\vector(0,1){30}}
\put(100,-5){\vector(-1,1){30}}
\put(100,25){\vector(-1,0){30}}

\put(62,-10){{\scriptsize A}}
\put(102,-10){{\scriptsize B}}
\put(62,27){{\scriptsize D}}
\put(102,27){{\scriptsize C}}

\put(170,-5){\circle*{3}}
\put(170,25){\circle*{3}}
\put(200,-5){\circle*{3}}
\put(200,25){\circle*{3}}
\put(170,-5){\vector(1,0){30}}
\put(170,-5){\vector(0,1){30}}
\put(200,-5){\vector(0,1){30}}
\put(200,25){\vector(-1,0){30}}

\put(240,-5){\circle*{3}}
\put(240,25){\circle*{3}}
\put(270,-5){\circle*{3}}
\put(270,25){\circle*{3}}
\put(240,-5){\vector(1,1){30}}
\put(240,-5){\vector(1,0){30}}
\put(270,-5){\vector(-1,1){30}}
\put(270,25){\vector(-1,0){30}}

\put(310,-5){\circle*{3}}
\put(310,25){\circle*{3}}
\put(340,-5){\circle*{3}}
\put(340,25){\circle*{3}}
\put(310,-5){\vector(1,1){30}}
\put(310,-5){\vector(0,1){30}}
\put(340,-5){\vector(0,1){30}}
\put(340,-5){\vector(-1,1){30}}

\put(180,-15){$-$}
\put(250,-15){$+$}
\put(320,-15){$+$}

\end{picture}
\end{center}
		\caption{An oriented chord diagram with labeled chords, its directed intersection graph,
and all the $4$-gons in it, with signs indicated}
		\label{feR2c}
\end{figure}
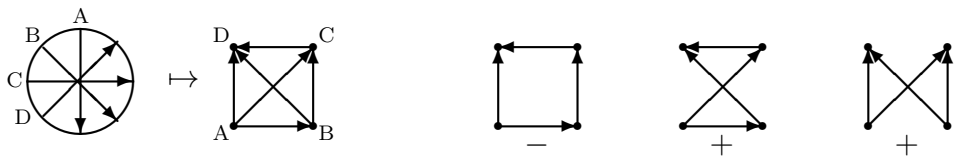

\begin{example}
Let us compute the value of~$R_2$ on the chord diagram~$d$ shown in Fig.~\ref{feR2c}.
Due to the Proposition above, an arbitrary orientation of the chords can be chosen.
For the one in the figure, the directed intersection graph~$\gamma(d)$ is shown in Fig.~\ref{feR2c}.
Figure~\ref{feR2c} also shows all the oriented $4$-gons in the intersection graph.
The difference between the number of positively and negatively oriented $4$-gons
is $2-1=1=R_2(d)$.
\end{example}

\subsection{Proof of the fact that~$R_k$ is a weight system}\label{sws}

\begin{theorem}
The function~$R_k$ on chord diagrams is indeed a weight system, that is, it satisfies
the $4$-term relation.
\end{theorem}

{\bf Proof.} The proof is similar to that of Proposition~\ref{pmod2}, but is slightly more complicated.
Let~$d$ be a chord diagram, and let~$A$ and~$B$ be a pair of chords in~$d$ having
neighboring ends and intersecting one another. Pick an orientation of the chords in~$d$ such that the arrow~$AB$
in the directed graph~$\gamma(d)$ of~$d$ is oriented from~$A$ to~$B$.

The difference $R_k(d)-R_k(d'_{AB})$ counts the number of signed circuits of length~$2k$
in~$\gamma(d)$ containing the edge~$AB$.

Each such circuit contains a sequence of vertices $CABD$. Let us split the set of such circuits into
the following three groups:
\begin{enumerate}
\item the vertex~$C$ is adjacent to~$B$ and~$D$ is adjacent to~$A$;
\item the vertex~$C$ is adjacent to~$B$, but~$D$ is not adjacent to~$A$;
\item the vertex C is not adjacent to B.
\end{enumerate}

All edge $2k$-gons in~$\widetilde{\gamma(d)}_{AB}$ passing through the four points $A,B,C,D$ admit a
similar classification.

Now, the $2k$-gons in~$\gamma(d)$ belonging to the second class are in one-to-one
correspondence with the $2k$-gons in~$\widetilde{\gamma(d)}_{AB}$ containing the path $CBAD$:
just preserve the other edges in each $2k$-gon.
In order to show that the signs of the corresponding circuits in both directed graphs coincide,
it suffices to consider all
possible mutual positions of the four chords $A,B,C,D$ in the chord diagram~$d$.
There are essentially two different such positions, depending on whether the chords~$C$ and~$D$
intersect one another. In both cases, for a chosen orientation of the chords,
we easily check coincidence of the signs of the corresponding $2k$-gons.


The
edge $2k$-gons of the third kind are the same in both graphs~$\gamma(d)$ and~$\widetilde\gamma(d)_{AB}$.

Finally, the edge $2k$-gons of the first kind in each of the two directed intersection graphs come in
pairs: the chain $CABD$ can be replaced with the chain $CBAD$. In order to prove that
the signs of the two $2k$-gons in a pair are opposite, it suffices to consider all
possible mutual positions of the four chords $A,B,C,D$ in the chord diagram~$d$.


Hence,
the required assertion follows.

\QED

\section{$R_k$ and intersection graphs}\label{sig}

Different chord diagrams can have one and the same intersection graph. Below,
we show that the value of the weight system~$R_k$ depends on the intersection
graph of the chord diagram rather than on the diagram itself. A natural
question then arises, whether~$R_k$ can be extended to a $4$-invariant of graphs.
We show that this is true for $k=2$ and $k=3$. The case of arbitrary~$k$ is discussed
in Sec.~\ref{s5}.

\subsection{Proof of the fact that~$R_k$ depends on the intersection graph only}\label{ssig}

\begin{theorem}\label{tig}
The weight system~$R_k$ acquires the same value on any two chord diagrams with
coinciding intersection graphs.
\end{theorem}

Since weight systems taking the same values on chord diagrams with
coinciding intersection graphs are in one-to one correspondence with finite
order knot invariants not distinguishing mutant knots~\cite{CL07}, we conclude the following:

\begin{corollary}
The canonical knot invariant associated to the weight system~$R_k$
does not distinguish mutant knots.
\end{corollary}

The proof of the theorem is based on a statement in~\cite{CL07} giving a complete
description of the situations where two chord diagrams have the same intersection graph.
We start with the definition of a share.

\begin{definition}
A \textit{share} is a part of a chord diagram consisting
of two arcs of the outer circle possessing the following property:
each chord one of whose ends belongs to these arcs
has both ends on these arcs.
\end{definition}

The complement of a share also is a share.
The whole chord diagram is its own share whose complement contains no chords.

\begin{figure}[ht]
\center{\includegraphics[width=.9\linewidth]{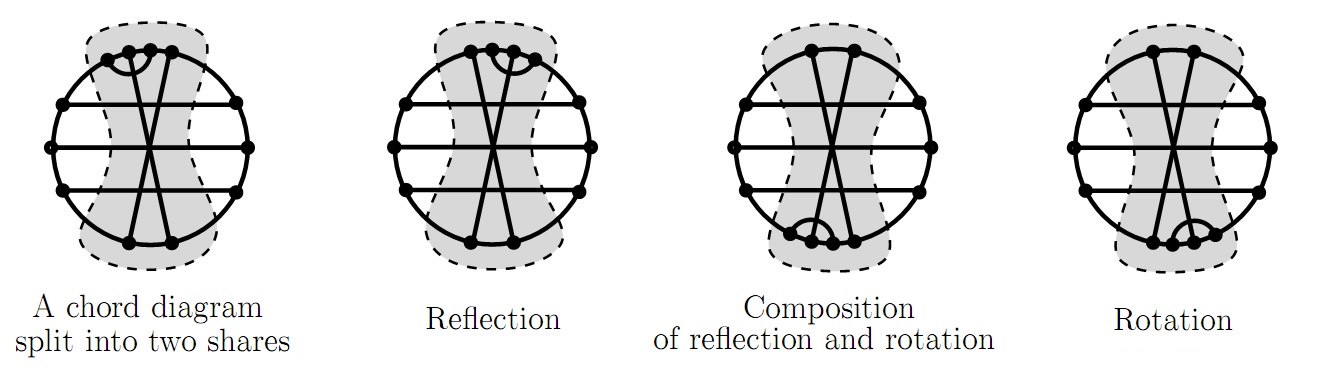}}
		\caption{Mutations of a chord diagram}
		\label{fm}
\end{figure}

\begin{definition}
A \textit{mutation of a chord diagram} is another chord diagram obtained by a rotation/reflection of a share;
see Fig.~\ref{fm}.
\end{definition}

Obviously, mutations preserve the intersection graphs of chord diagrams.
We call the subgraph of the intersection graph induced by the vertices corresponding
to chords forming a given share also a share in the intersection graph.

\begin{theorem}{\rm\cite{CL07}} \label{cd-mut-theorem}
Two chord diagrams have the same intersection graph if and only if
they are related by a sequence of mutations.
\end{theorem}

{\bf Proof} of Theorem~\ref{tig}. We are going to prove that a mutation of a
chord diagram does not change the value of the weight system~$R_k$,
for arbitrary $k=2,3,4,\dots$. It suffices to consider only reflection and rotation
of a share; for the composition of reflection and rotation, the result will follow
automatically.

Pick an arbitrary orientation of chords in a chord diagram.
A mutation of a chord diagram with oriented chords produces a chord diagram with oriented chords
in a natural way.
For a given share in the intersection graph, its reflection results in inverting the orientation of all the arrows
belonging to the corresponding share in the intersection graph.
Similarly, rotation of a share results in inverting all the arrows
between the two shares in the intersection graph.


Denote the two sets of vertices belonging to the complementary shares in the intersection graph
by~$U$ and~$W$, respectively. The set~$U$ contains a subset~$u\subset U$ and the set $W$
contains a subset~$w\subset W$ such that
\begin{itemize}
\item any vertex from~$u$ is connected to all vertices of~$w$ (and vice versa);
\item each edge connecting vertices from~$U$ and~$W$ connects, in fact, a vertex from~$u$
to a vertex from~$w$.
\end{itemize}
Denote the set of all edges in the intersection graph connecting a vertex from~$u$
to a vertex from~$w$ by~$K(u,w)$. The graph with the set of vertices $u\sqcup w$
and the set of edges $K(u,w)$ is the complete bipartite graph with the parts $u,w$.

Any circuit in the intersection graph contains an even number of edges from~$K(u,w)$.
Indeed, any path starting in~$U$ switches between~$U$ and~$W$ after each passing through
an edge in~$K(u,w)$. Since a circuit returns to the original vertex, the number
of such passings must be even.

Now, the rotation mutation changes the orientation of all arrows in~$K(u,w)$, whence
of an even number of arrows in any circuit. In particular, it changes orientation
of an even number of arrows in any $2k$-gon, hence preserving its sign.

\begin{figure}[ht]
\center{\includegraphics[width=.9\linewidth]{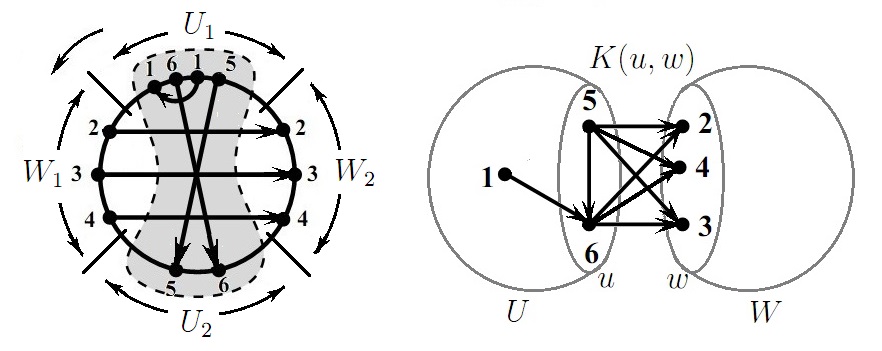}}
		\caption{Distinguished orientation of chords allowing for a
specific orientation of share connecting arrows in the directed intersection graph}
		\label{forient}
\end{figure}

Pick an arbitrary orientation of the chord diagram such that all the arrows in~$K(u,w)$ are oriented
from~$U$ to~$W$. Such an orientation always exists. Indeed, if the share~$U$ consists of two arcs~$U_1,U_2$
and the share~$W$ consists of two arcs~$W_1,W_2$ following along the positive direction of the outer circle
in the alternating order, $U_1,W_1,U_2,W_2$, see Fig.~\ref{forient}, then it suffices to orient
all the chords connecting~$U_1,U_2$ and $W_1,W_2$ from~$U_1$ to~$U_2$ and from~$W_1$ to~$W_2$,
respectively.

What happens if we reflect the share~$W$? If the $2k$-gon we consider has even number
of arrows in the share~$W$, then its sign remains the same.
The circuits with odd number of arrows in~$W$ split into pairs in the following way.
The circuit intersects the set $K(u,w)$ by even number of arrows. Let us
split these arrows into pairs: two edges belong to one and the same pair
iff their ends in the share~$W$ are connected by a path that is a part of the circuit totally
contained in~$W$. For such pair of arrows, let the first one be $u_1w_1$
and the second one be $u_2w_2$, $u_i\in U$, $w_i\in W$, $i=1,2$ (it may well happen
that either $w_1=w_2$ or $u_1=u_2$, but not both).

\begin{figure}[ht]
\center{\includegraphics[width=.5\linewidth]{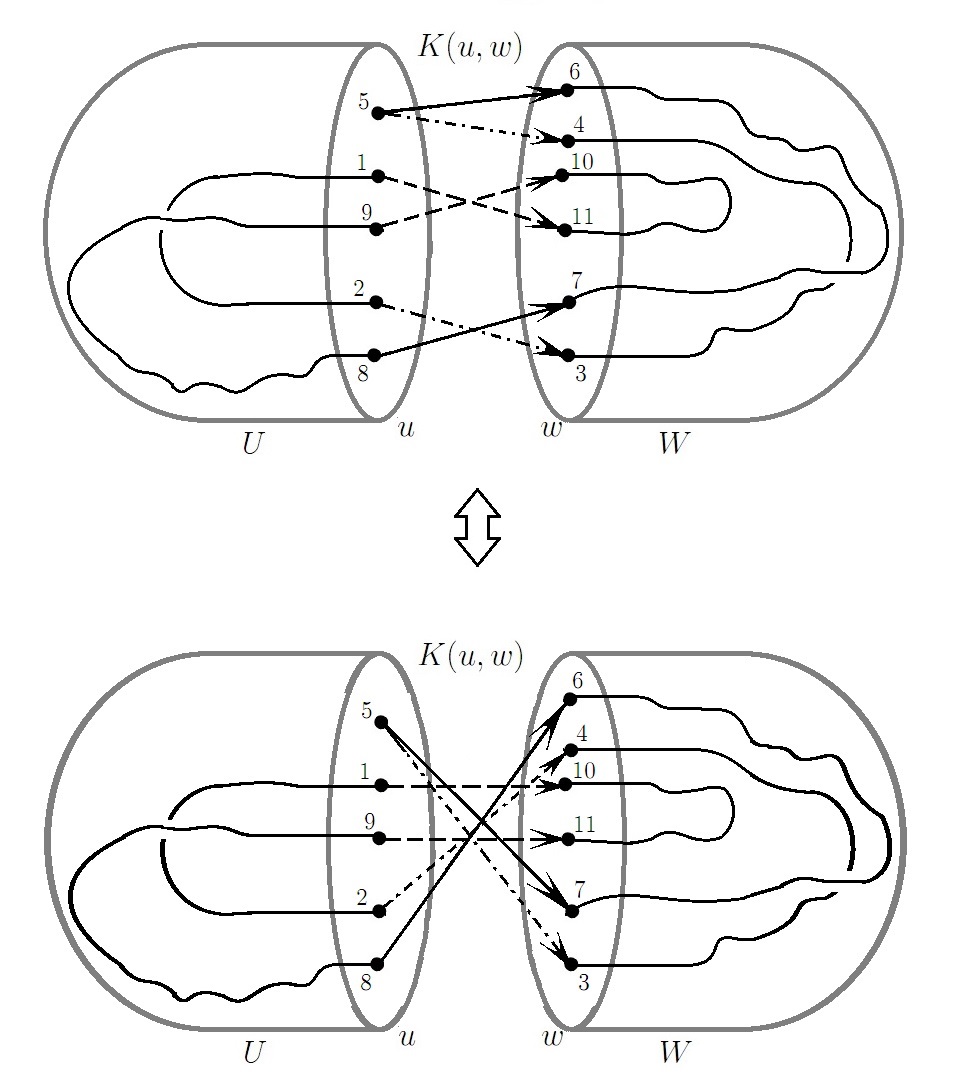}}
		\caption{Two circuits obtained from one another by switching pairs}
		\label{fsp}
\end{figure}

By {\it switching pairs\/} we mean replacing each pair of edges $u_1w_1$,
$u_2w_2$, with the pair $u_1w_2$, $u_2w_1$; see Fig.~\ref{fsp}.

Switching all pairs for a circuit produces a new circuit. This transformation is
an involution: switching all pairs in this new circuit restores the original one.
Thus, all circuits are split into pairs. We are going to show that
any two circuits of even length~$2k$ having odd number of arrows in~$W$ and
belonging to the same pair have opposite signs. Indeed, if we choose directions
in both circuits that coincide inside the share~$U$, then they are opposite inside
the share~$W$. Both circuits have the same number of edges inside~$K(u,w)$
in either direction. Therefore, the only difference in the orientation is
inside the share~$W$, where each arrow has opposite directions in the two circuits.
Since the number of arrows in~$W$ is odd, the signs
of the two circuits are opposite.

Thus, both kinds of mutation preserve the value of~$R_k$. \QED

\subsection{Extending~$R_k$ to arbitrary graphs}\label{ss4i}
Theorem~\ref{tig} above shows that each weight system~$R_k$ defines a function
on intersection graphs; we denote this function by~$R_k$ as well.
In this section we prove the following special case of a general theorem proved in Sec.~\ref{s5}.

\begin{proposition}\label{pge}
The functions~$R_2$ and~$R_3$ can be extended to $4$-invariants of graphs.
\end{proposition}

{\bf Proof.} In order to prove that~$R_k$ can be extended to a $4$-invariant
of graphs, it suffices to prove this for graphs with exactly $2k$ vertices.
Indeed, denote by $\overline{R}_k$ the function whose values
on graphs with~$2k$ vertices coincides with the constructed extension, and
which vanishes on all other graphs. Then the convolution $\overline{R}_k * U$
provides us with the desired extension to arbitrary graphs; here~$U$
is the function whose value on the arbitrary graph is equal to~$1$.
It is obvious that when restricted to intersection graphs,  $\overline{R}_k * U$
coincides with~$R_k$. Here we make use of the fact that the vector space
spanned by graphs carries a natural
graded bialgebra structure: the product of two graphs is given by their disjoint union,
and the coproduct $\mu(\G)$ of a graph~$\G$ is given by
$$
\mu(\G)=\sum_{I\subset V(\G)}\G|_I\otimes \G|_{V(\G)\setminus I},
$$
where~$\G|_I$ is the subgraph of~$\G$ induced by a subset~$I$ of its vertices.
The convolution product, denoted by~$*$, is induced on the
dual space of graph invariants from the coproduct on the space of graphs~\cite{L00}:
$$
f*g(\G)=f\otimes g(\mu(G)),
$$
for arbitrary graph invariants~$f$ and~$g$. These operations naturally descend
to the $4$-bialgebra of graphs and $4$-invariants, respectively.

If~$k=2$, then each graph with~$2k=4$ vertices is an intersection graph, and
any $4$-term relation for intersection graphs has a chord diagram counterpart,
meaning we are done.

\begin{figure}[ht]

\begin{center}
\begin{picture}(300,450)(10,20)

\thicklines

\put(55,551){\circle*{4}}
\put(55,575){\circle*{4}}
\put(34,562){\circle*{4}}
\put(76,562){\circle*{4}}
\put(40,535){\circle*{4}}
\put(70,535){\circle*{4}}
\put(55,550){\line(0,1){24}}
\put(55,551){\line(2,1){21}}
\put(55,551){\line(-2,1){21}}
\put(55,550){\line(1,-1){16}}
\put(55,550){\line(-1,-1){16}}
\put(40,535){\line(1,0){32}}
\put(55,575){\line(2,-1){21}}
\put(55,575){\line(-2,-1){21}}
\put(40,535){\line(-1,4){7}}
\put(70,535){\line(1,4){7}}

\put(90,551){$-$}
\put(28,532){$A$}
\put(72,532){$B$}

\put(135,551){\circle*{4}}
\put(135,575){\circle*{4}}
\put(114,562){\circle*{4}}
\put(156,562){\circle*{4}}
\put(120,535){\circle*{4}}
\put(150,535){\circle*{4}}
\put(135,550){\line(0,1){24}}
\put(135,551){\line(2,1){21}}
\put(135,551){\line(-2,1){21}}
\put(135,550){\line(1,-1){16}}
\put(135,550){\line(-1,-1){16}}
\put(135,575){\line(2,-1){21}}
\put(135,575){\line(-2,-1){21}}
\put(120,535){\line(-1,4){7}}
\put(150,535){\line(1,4){7}}

\put(170,551){$=$}
\put(108,532){$A$}
\put(152,532){$B$}

\put(215,551){\circle*{4}}
\put(215,575){\circle*{4}}
\put(194,562){\circle*{4}}
\put(236,562){\circle*{4}}
\put(200,535){\circle*{4}}
\put(230,535){\circle*{4}}
\put(215,550){\line(0,1){24}}
\put(215,551){\line(2,1){21}}
\put(215,551){\line(-2,1){21}}
\put(215,550){\line(1,-1){16}}
\put(200,536){\line(3,2){36}}
\put(200,535){\line(1,0){32}}
\put(215,575){\line(2,-1){21}}
\put(215,575){\line(-2,-1){21}}
\put(200,535){\line(-1,4){7}}
\put(230,535){\line(1,4){7}}

\put(250,551){$-$}
\put(188,532){$A$}
\put(232,532){$B$}

\put(295,551){\circle*{4}}
\put(295,575){\circle*{4}}
\put(274,562){\circle*{4}}
\put(316,562){\circle*{4}}
\put(280,535){\circle*{4}}
\put(310,535){\circle*{4}}
\put(295,550){\line(0,1){24}}
\put(295,551){\line(2,1){21}}
\put(295,551){\line(-2,1){21}}
\put(295,550){\line(1,-1){16}}
\put(280,536){\line(3,2){36}}
\put(295,575){\line(2,-1){21}}
\put(295,575){\line(-2,-1){21}}
\put(280,535){\line(-1,4){7}}
\put(310,535){\line(1,4){7}}

\put(268,532){$A$}
\put(312,532){$B$}

\put(128,515){$-1$}
\put(205,515){$-3$}
\put(285,515){$-1$}

\put(55,451){\circle*{4}}
\put(55,475){\circle*{4}}
\put(34,462){\circle*{4}}
\put(76,462){\circle*{4}}
\put(40,435){\circle*{4}}
\put(70,435){\circle*{4}}
\put(55,450){\line(0,1){24}}
\put(55,451){\line(2,1){21}}
\put(55,451){\line(-2,1){21}}
\put(55,450){\line(1,-1){16}}
\put(55,450){\line(-1,-1){16}}
\put(40,435){\line(1,0){32}}
\put(55,475){\line(2,-1){21}}
\put(55,475){\line(-2,-1){21}}
\put(40,435){\line(-1,4){7}}
\put(70,435){\line(1,4){7}}

\put(90,451){$-$}
\put(28,432){$A$}
\put(55,457){$B$}

\put(135,451){\circle*{4}}
\put(135,475){\circle*{4}}
\put(114,462){\circle*{4}}
\put(156,462){\circle*{4}}
\put(120,435){\circle*{4}}
\put(150,435){\circle*{4}}
\put(135,450){\line(0,1){24}}
\put(135,451){\line(2,1){21}}
\put(135,451){\line(-2,1){21}}
\put(135,450){\line(1,-1){16}}
\put(120,435){\line(1,0){32}}
\put(135,475){\line(2,-1){21}}
\put(135,475){\line(-2,-1){21}}
\put(120,435){\line(-1,4){7}}
\put(150,435){\line(1,4){7}}

\put(170,451){$=$}
\put(108,432){$A$}
\put(135,457){$B$}

\put(215,451){\circle*{4}}
\put(215,475){\circle*{4}}
\put(194,462){\circle*{4}}
\put(236,462){\circle*{4}}
\put(200,435){\circle*{4}}
\put(230,435){\circle*{4}}
\put(215,450){\line(0,1){24}}
\put(215,451){\line(2,1){21}}
\put(215,451){\line(-2,1){21}}
\put(215,450){\line(1,-1){16}}
\put(215,450){\line(-1,-1){16}}
\put(201,436){\line(1,3){13}}
\put(200,436){\line(3,2){36}}
\put(215,475){\line(2,-1){21}}
\put(215,475){\line(-2,-1){21}}
\put(230,435){\line(1,4){7}}

\put(250,451){$-$}
\put(188,432){$A$}
\put(215,457){$B$}

\put(295,451){\circle*{4}}
\put(295,475){\circle*{4}}
\put(274,462){\circle*{4}}
\put(316,462){\circle*{4}}
\put(280,435){\circle*{4}}
\put(310,435){\circle*{4}}
\put(295,450){\line(0,1){24}}
\put(295,451){\line(2,1){21}}
\put(295,451){\line(-2,1){21}}
\put(295,450){\line(1,-1){16}}
\put(280,436){\line(3,2){36}}
\put(281,436){\line(1,3){13}}
\put(295,475){\line(2,-1){21}}
\put(295,475){\line(-2,-1){21}}
\put(310,435){\line(1,4){7}}

\put(268,432){$A$}
\put(295,457){$B$}

\put(128,415){$-1$}
\put(205,415){$-1$}
\put(292,415){$1$}

\put(55,351){\circle*{4}}
\put(55,375){\circle*{4}}
\put(34,362){\circle*{4}}
\put(76,362){\circle*{4}}
\put(40,335){\circle*{4}}
\put(70,335){\circle*{4}}
\put(55,350){\line(0,1){24}}
\put(55,351){\line(2,1){21}}
\put(55,351){\line(-2,1){21}}
\put(55,350){\line(1,-1){16}}
\put(55,350){\line(-1,-1){16}}
\put(40,335){\line(1,0){32}}
\put(55,375){\line(2,-1){21}}
\put(55,375){\line(-2,-1){21}}
\put(40,335){\line(-1,4){7}}
\put(70,335){\line(1,4){7}}

\put(90,351){$-$}
\put(28,332){$B$}
\put(55,357){$A$}

\put(135,351){\circle*{4}}
\put(135,375){\circle*{4}}
\put(114,362){\circle*{4}}
\put(156,362){\circle*{4}}
\put(120,335){\circle*{4}}
\put(150,335){\circle*{4}}
\put(135,350){\line(0,1){24}}
\put(135,351){\line(2,1){21}}
\put(135,351){\line(-2,1){21}}
\put(135,350){\line(1,-1){16}}
\put(120,335){\line(1,0){32}}
\put(135,375){\line(2,-1){21}}
\put(135,375){\line(-2,-1){21}}
\put(120,335){\line(-1,4){7}}
\put(150,335){\line(1,4){7}}

\put(170,351){$=$}
\put(108,332){$B$}
\put(135,357){$A$}

\put(215,351){\circle*{4}}
\put(215,375){\circle*{4}}
\put(194,362){\circle*{4}}
\put(236,362){\circle*{4}}
\put(200,335){\circle*{4}}
\put(230,335){\circle*{4}}
\put(215,350){\line(0,1){24}}
\put(215,351){\line(2,1){21}}
\put(215,350){\line(-1,-1){16}}
\put(200,335){\line(1,0){32}}
\put(215,375){\line(2,-1){21}}
\put(215,375){\line(-2,-1){21}}
\put(200,335){\line(-1,4){7}}
\put(230,335){\line(1,4){7}}

\put(250,351){$-$}
\put(188,332){$B$}
\put(215,357){$A$}

\put(295,351){\circle*{4}}
\put(295,375){\circle*{4}}
\put(274,362){\circle*{4}}
\put(316,362){\circle*{4}}
\put(280,335){\circle*{4}}
\put(310,335){\circle*{4}}
\put(295,350){\line(0,1){24}}
\put(295,351){\line(2,1){21}}
\put(280,335){\line(1,0){32}}
\put(295,375){\line(2,-1){21}}
\put(295,375){\line(-2,-1){21}}
\put(280,335){\line(-1,4){7}}
\put(310,335){\line(1,4){7}}

\put(268,332){$B$}
\put(295,357){$A$}

\put(128,315){$-1$}
\put(205,315){$-1$}
\put(292,315){$1$}

\put(55,270){\circle*{4}}
\put(55,240){\circle*{4}}
\put(35,280){\circle*{4}}
\put(75,280){\circle*{4}}
\put(35,230){\circle*{4}}
\put(75,230){\circle*{4}}
\put(35,230){\line(0,1){50}}
\put(35,230){\line(1,0){40}}
\put(35,280){\line(1,0){40}}
\put(75,230){\line(0,1){50}}
\put(55,240){\line(0,1){30}}
\put(35,230){\line(2,1){20}}
\put(35,280){\line(2,-1){20}}
\put(75,280){\line(-2,-1){20}}
\put(75,230){\line(-2,1){20}}

\put(90,251){$-$}
\put(43,240){$A$}
\put(58,265){$B$}

\put(135,270){\circle*{4}}
\put(135,240){\circle*{4}}
\put(115,280){\circle*{4}}
\put(155,280){\circle*{4}}
\put(115,230){\circle*{4}}
\put(155,230){\circle*{4}}
\put(115,230){\line(0,1){50}}
\put(115,230){\line(1,0){40}}
\put(115,280){\line(1,0){40}}
\put(155,230){\line(0,1){50}}
\put(115,230){\line(2,1){20}}
\put(115,280){\line(2,-1){20}}
\put(155,280){\line(-2,-1){20}}
\put(155,230){\line(-2,1){20}}

\put(170,251){$=$}
\put(123,240){$A$}
\put(138,265){$B$}

\put(215,270){\circle*{4}}
\put(215,240){\circle*{4}}
\put(195,280){\circle*{4}}
\put(235,280){\circle*{4}}
\put(195,230){\circle*{4}}
\put(235,230){\circle*{4}}
\put(195,230){\line(0,1){50}}
\put(195,230){\line(1,0){40}}
\put(195,280){\line(1,0){40}}
\put(235,230){\line(0,1){50}}
\put(215,240){\line(0,1){30}}
\put(195,230){\line(2,1){20}}
\put(195,280){\line(2,-1){20}}
\put(235,280){\line(-2,-1){20}}
\put(235,230){\line(-2,1){20}}
\put(195,280){\line(1,-2){20}}
\put(235,280){\line(-1,-2){20}}

\put(250,251){$-$}
\put(203,240){$A$}
\put(218,265){$B$}

\put(295,270){\circle*{4}}
\put(295,240){\circle*{4}}
\put(275,280){\circle*{4}}
\put(315,280){\circle*{4}}
\put(275,230){\circle*{4}}
\put(315,230){\circle*{4}}
\put(275,230){\line(0,1){50}}
\put(275,230){\line(1,0){40}}
\put(275,280){\line(1,0){40}}
\put(315,230){\line(0,1){50}}
\put(275,230){\line(2,1){20}}
\put(275,280){\line(2,-1){20}}
\put(315,280){\line(-2,-1){20}}
\put(315,230){\line(-2,1){20}}
\put(275,280){\line(1,-2){20}}
\put(315,280){\line(-1,-2){20}}

\put(283,240){$A$}
\put(298,265){$B$}

\put(132,215){$1$}
\put(205,215){$-3$}
\put(285,215){$-1$}

\put(55,170){\circle*{4}}
\put(55,140){\circle*{4}}
\put(35,180){\circle*{4}}
\put(75,180){\circle*{4}}
\put(35,130){\circle*{4}}
\put(75,130){\circle*{4}}
\put(35,130){\line(0,1){50}}
\put(35,130){\line(1,0){40}}
\put(35,180){\line(1,0){40}}
\put(75,130){\line(0,1){50}}
\put(55,140){\line(0,1){30}}
\put(35,130){\line(2,1){20}}
\put(35,180){\line(2,-1){20}}
\put(75,180){\line(-2,-1){20}}
\put(75,130){\line(-2,1){20}}

\put(90,151){$-$}
\put(43,140){$A$}
\put(78,125){$B$}

\put(135,170){\circle*{4}}
\put(135,140){\circle*{4}}
\put(115,180){\circle*{4}}
\put(155,180){\circle*{4}}
\put(115,130){\circle*{4}}
\put(155,130){\circle*{4}}
\put(115,130){\line(0,1){50}}
\put(115,130){\line(1,0){40}}
\put(115,180){\line(1,0){40}}
\put(155,130){\line(0,1){50}}
\put(135,140){\line(0,1){30}}
\put(115,130){\line(2,1){20}}
\put(115,180){\line(2,-1){20}}
\put(155,180){\line(-2,-1){20}}

\put(170,151){$=$}
\put(123,140){$A$}
\put(158,125){$B$}

\put(215,170){\circle*{4}}
\put(215,140){\circle*{4}}
\put(195,180){\circle*{4}}
\put(235,180){\circle*{4}}
\put(195,130){\circle*{4}}
\put(235,130){\circle*{4}}
\put(195,130){\line(0,1){50}}
\put(195,130){\line(1,0){40}}
\put(195,180){\line(1,0){40}}
\put(235,130){\line(0,1){50}}
\put(215,140){\line(0,1){30}}
\put(195,180){\line(2,-1){20}}
\put(235,180){\line(-2,-1){20}}
\put(235,130){\line(-2,1){20}}
\put(235,180){\line(-1,-2){20}}
\put(235,180){\line(-1,-2){20}}

\put(250,151){$-$}
\put(203,140){$A$}
\put(238,125){$B$}

\put(295,170){\circle*{4}}
\put(295,140){\circle*{4}}
\put(275,180){\circle*{4}}
\put(315,180){\circle*{4}}
\put(275,130){\circle*{4}}
\put(315,130){\circle*{4}}
\put(275,130){\line(0,1){50}}
\put(275,130){\line(1,0){40}}
\put(275,180){\line(1,0){40}}
\put(315,130){\line(0,1){50}}
\put(295,140){\line(0,1){30}}
\put(275,180){\line(2,-1){20}}
\put(315,180){\line(-2,-1){20}}
\put(315,180){\line(-1,-2){20}}

\put(283,140){$A$}
\put(318,125){$B$}

\put(130,115){$-1$}
\put(205,115){$-1$}
\put(285,115){$-1$}

\end{picture}
\end{center}

		\caption{$4$-term relations for the $5$-wheel and the prism expressing them in terms
of intersection graphs; the value of the weight system~$R_3$ is indicated}
		\label{fR3}
\end{figure}
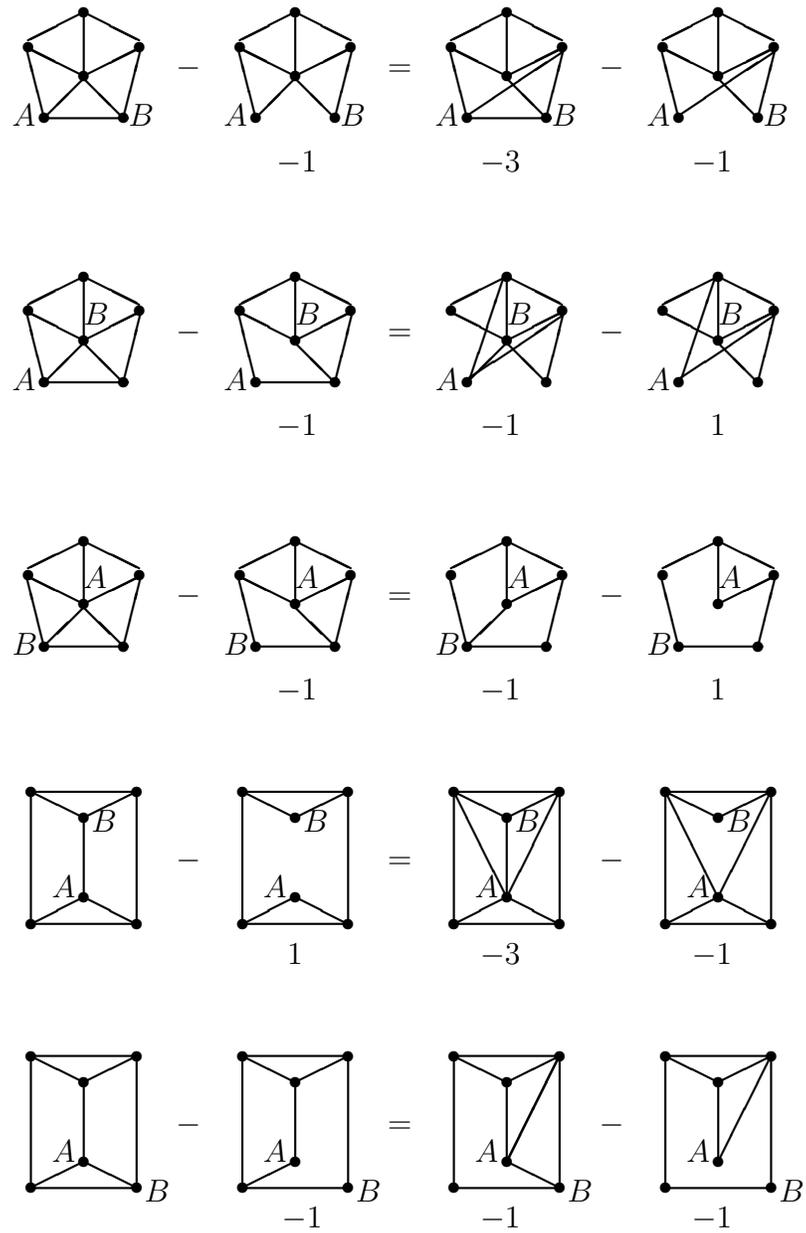

For $k=3$, there are two graphs with 6 vertices that are not intersection graphs,
namely, the $5$-wheel and the $3$-prism; see Fig.~\ref{fR3}. There are three ways to express
the $5$-wheel as a linear combination of intersection graphs through the $4$-term
relation, and for the $3$-prism there are two such ways; see Fig.~\ref{fR3}.
A direct computation shows that in both cases all the different representations
provide the same value, which one can admit for the value of the extended graph
invariant~$R_3$: for the $5$-wheel it is~$-3$, and for the $3$-prism it is~$-1$.
Now, a direct verification, taking into account those $4$-term relations for graphs
that do not have chord diagram counterparts, ensures that this extended invariant of graphs
with~$6$ vertices is indeed a $4$-invariant. \QED

\section{$R_k$ and the $\fsl_2$-weight system}\label{ssl2}

It is well known that for an arbitrary Lie algebra~$\g$ endowed with a
nondegenerate invariant scalar product $(\cdot,\cdot)$ one can construct a weight system~$b^{\g}$
with values in the center $ZU(\g)$ of the universal enveloping algebra $U(\g)$.
Invariancy of the scalar product means that $(x,[y,z])=([x,y],z)$ for any
three elements $x,y,z\in\g$.

The $3$-dimensional Lie algebra~$\fsl_2$ (taken, for definiteness, over the field~$\C$ of complex numbers)
provides the first nontrivial example of this construction. Since the center $ZU(\fsl_2)$
of the universal enveloping algebra $U(\fsl_2)$ is isomorphic to the ring of polynomials
in the Casimir element~$c$ of~$\fsl_2$, we obtain a weight system with values
in the ring $\C[c]$ of polynomials in a single variable~$c$; see \cite{BN95} and \cite{K93}. We denote this weight system
by $b^{\fsl_2}:\cM\to\C[c]$.
Here~$\cM$ denotes the graded vector space spanned by chord diagrams modulo the $4$-term relations;
$$
\cM=\cM_0\oplus\cM_1\oplus\cM_2\oplus\dots,
$$
where the finite dimensional vector space~$\cM_n$ is spanned by chord diagrams with~$n$
chords, modulo the $4$-term relations. The space~$\cM$ is a graded commutative cocommutative
Hopf algebra, with a multiplication $m:\cM\otimes\cM\to\cM$, and a comultiplication $\mu:\cM\to\cM\otimes\cM$.
The weight system $b^{\fsl_2}$ is an algebra homomorphism.
It takes an arbitrary chord diagram
with~$n$ chords to a polynomial of degree~$n$ in~$c$ with the leading coefficient~$1$.

The invariant nondegenerate scalar product~$(\cdot,\cdot)$ on~$\fsl_2$ can be chosen in a unique
way up to a nonzero multiplicative constant. The choice of the multiplicative
constant affects the non-leading coefficients of the value of~$b^{\fsl_2}$ on a chord
diagram. Below, we use the normalization chosen
in~\cite{LZ04}, which differs from that in~\cite{CV97} and~\cite{CDM12}.

Being a graded commutative cocommutative Hopf algebra, $\cM$ is generated by its primitive
elements. Recall that an element~$p$ of a Hopf algebra is said to be {\it primitive\/}
if  $\mu(p)=1\otimes p+p\otimes1$. Primitive elements form a graded vector subspace~$\cPM$ of~$\cM$,
$$
\cPM=\cPM_1\oplus\cPM_2\oplus\cPM_3\oplus\dots,\qquad \cPM_n\subset\cM_n.
$$
For each $n=1,2,3,\dots$, there exists a natural projection $\pi_n:\cM_n\to\cPM_n$
of the space of chord diagrams of order~$n$ to the subspace~$\cPM_n$ of primitive elements
along the subspace $\cDM_n$ spanned by decomposable chord diagrams. A chord diagram with~$n$
chords is said to be {\it decomposable\/} if it can be represented as a product of
two chord diagrams, each with less than~$n$ chords.

The following statement is an immediate consequence of a result obtained in~\cite{CV97}, Theorem~3,
for an arbitrary simple Lie algebra.

\begin{proposition}~{\rm \cite{CV97}}
The value of the weight system $b^{\fsl_2}$ on an arbitrary element in~$\cPM_n$
is a polynomial of degree at most $n/2$ in~$c$.
\end{proposition}

Now we are able to formulate the conjecture mentioned in the introduction, which
relates the weight systems~$R_k$ and $b^{\fsl_2}$.

\begin{conjecture}\label{conj-sl2-hc}
For any chord diagram~$d$ with~$2k$ chords, the coefficient of~$c^k$ in the value
$b^{\fsl_2}(\pi_{2k}(d))$ is $2R_k(d)$.
\end{conjecture}

We did not manage to prove the Conjecture, but the computer experiments show that
it is true up to $k=4$. Below, we present Table 1, which contains the values
of the weight systems $R_k(\cdot)$ and $b^{sl_2}(\pi_{2k}(\cdot))$ for some chord
diagrams with~$2k$ chords.

The original definition of a weight system~$b^\g$ is too cumbersome from
the computational point of view, which makes it difficult to be computed
even for chord diagrams with very few chords.
The computations we made used the recurrence relation for~$b^{\fsl_2}$
obtained in~\cite{CV97}, Theorem~2:

\begin{proposition}
If a chord diagram contains a leaf, that is, a chord intersecting only
one other chord, then the value of the $\fsl_2$ universal weight
system on the diagram is~$(c-1)$ times its value on the result of
deleting the leaf. In addition, the value of the weight system $b^{\fsl_2}$ on
a chord diagram satisfies the recurrence relations shown in Fig.~\ref{7fig:31}.
\end{proposition}

\begin{figure}[htbp]
%
%
\def\bepi#1{\makebox[33pt]{\unitlength=18pt
            \begin{picture}(1.8,1.1)(-0.98,-0.2)  #1
            \end{picture}} }
\def\sctw#1#2#3#4{
   \bezier{25}(-0.26,0.97)(0,1.035)(0.26,0.97)
   \bezier{4}(0.26,0.97)(0.52,0.9)(0.71,0.71)
   \bezier{#1}(0.71,0.71)(0.9,0.52)(0.97,0.26)
   \bezier{4}(0.97,0.26)(1.035,0)(0.97,-0.26)
   \bezier{#2}(0.97,-0.26)(0.9,-0.52)(0.71,-0.71)
   \bezier{4}(0.71,-0.71)(0.52,-0.9)(0.26,-0.97)
   \bezier{25}(0.26,-0.97)(0,-1.035)(-0.26,-0.97)
   \bezier{4}(-0.26,-0.97)(-0.52,-0.9)(-0.71,-0.71)
   \bezier{#3}(-0.71,-0.71)(-0.9,-0.52)(-0.97,-0.26)
   \bezier{4}(-0.97,-0.26)(-1.035,0)(-0.97,0.26)
   \bezier{#4}(-0.97,0.26)(-0.9,0.52)(-0.71,0.71)
   \bezier{4}(-0.71,0.71)(-0.52,0.9)(-0.26,0.97)
}
%
\def\dslone{\bepi{\sctw{25}{25}{4}{4} \slchone \slchtwo \slchthree}}
\def\dsltwo{\bepi{\sctw{25}{25}{4}{4} \slchone \slchfour \slchthree}}
\def\dslthree{\bepi{\sctw{25}{25}{4}{4} \slchone \slchtwo \slchfive}}
\def\dslfour{\bepi{\sctw{25}{25}{4}{4} \slchone \slchfour \slchfive}}
\def\dslfive{\bepi{\sctw{25}{25}{4}{4} \slchsix \slchtth}}
\def\dslsix{\bepi{\sctw{25}{25}{4}{4} \slchseven \slcheight}}
\def\dslseven{\bepi{\sctw{4}{25}{4}{25} \slchone \slchnine \slchthree}}
\def\dsleight{\bepi{\sctw{4}{25}{4}{25} \slchone \slchten \slchthree}}
\def\dslnine{\bepi{\sctw{4}{25}{4}{25} \slchone \slchnine \slchfive}}
\def\dslten{\bepi{\sctw{4}{25}{4}{25} \slchone \slchten \slchfive}}
\def\dslel{\bepi{\sctw{4}{25}{4}{25} \slchel \put(-0.866,0.5){\circle*{0.15}}
    \put(0.866,-0.5){\circle*{0.15}} \put(0.866,-0.5){\line(-5,3){1.72}} }}
\def\dsltw{\bepi{\sctw{4}{25}{4}{25} \slchtw \slcheight}}
%
%
\def\aslone{\bepi{\sctw{25}{25}{4}{4} \slchone \slchseven \slcheight}}
\def\asltwo{\bepi{\sctw{25}{25}{4}{4} \slchone \slchseven \achone}}
\def\aslthree{\bepi{\sctw{25}{25}{4}{4} \slchone \achtwo \slcheight}}
\def\aslfour{\bepi{\sctw{25}{25}{4}{4} \slchone \achone \achtwo}}
\def\aslfive{\bepi{\sctw{25}{25}{4}{4} \slchtwo \slchthree}}
\def\asltse{\bepi{\sctw{25}{4}{25}{4} \slchone \achsix \slchtwo}}
\def\asltei{\bepi{\sctw{25}{4}{25}{4} \slchone \achsix \slchfour}}
\def\asltni{\bepi{\sctw{25}{4}{25}{4} \slchone \achfive \slchtwo}}
\def\asltte{\bepi{\sctw{25}{4}{25}{4} \slchone \achfive \slchfour}}
\def\asltto{\bepi{\sctw{25}{4}{25}{4} \achnine \slchseven}}
\def\aslttt{\bepi{\sctw{25}{4}{25}{4}\slchel\put(-0.866,-0.5){\circle*{0.15}}
   \put(0.866,0.5){\circle*{0.15}}\put(0.866,0.5){\line(-5,-3){1.72}} }}
%
%
\def\di#1{\raisebox{0pt}[20pt][23pt]{#1}}
\def\wdi#1{W\biggl( {\raisebox{0pt}[20pt][30pt]{#1}} \biggr)}
\def\spwdi#1{W\biggl( \mbox{#1} \biggr)}
\def\slrel#1#2#3#4#5#6#7{
\begin{eqnarray*}
\wdi{#1}-\wdi{#2}-\wdi{#3}+\wdi{#4} = \makebox[30pt]{}&  \\
\hfill = 2 \wdi{#5} - 2 \wdi{#6} #7 &
\end{eqnarray*} }
%
%
%
\def\slchone{\put(0,1){\circle*{0.15}}\put(0,-1){\circle*{0.15}}
             \put(0,-1){\line(0,1){2}} }
\def\slchtwo{\put(-0.174,0.985){\circle*{0.15}}
             \put(0.866,0.5){\circle*{0.15}}
             \bezier{70}(-0.174,0.985)(0.21,0.45)(0.866,0.5) }
\def\slchthree{\put(-0.174,-0.985){\circle*{0.15}}
               \put(0.866,-0.5){\circle*{0.15}}
               \bezier{70}(-0.174,-0.985)(0.21,-0.45)(0.866,-0.5) }
\def\slchfour{\put(0.174,0.985){\circle*{0.15}}
              \put(0.866,0.5){\circle*{0.15}}
              \bezier{60}(0.174,0.985)(0.36,0.48)(0.866,0.5) }
\def\slchfive{\put(0.174,-0.985){\circle*{0.15}}
              \put(0.866,-0.5){\circle*{0.15}}
              \bezier{60}(0.174,-0.985)(0.36,-0.48)(0.866,-0.5) }
\def\slchsix{\put(-0.174,0.985){\circle*{0.15}}
             \put(-0.174,-0.985){\circle*{0.15}}
             \bezier{100}(-0.174,0.985)(0.3,0)(-0.174,-0.985) }
\def\slchseven{\put(-0.174,-0.985){\circle*{0.15}}
               \put(0.866,0.5){\circle*{0.15}}
               \bezier{90}(-0.174,-0.985)(0,0)(0.866,0.5) }
\def\slcheight{\put(-0.174,0.985){\circle*{0.15}}
               \put(0.866,-0.5){\circle*{0.15}}
               \bezier{90}(-0.174,0.985)(0,0)(0.866,-0.5) }
\def\slchnine{\put(-0.174,0.985){\circle*{0.15}}
              \put(-0.866,0.5){\circle*{0.15}}
              \bezier{60}(-0.174,0.985)(-0.36,0.48)(-0.866,0.5) }
\def\slchten{\put(0.174,0.985){\circle*{0.15}}
             \put(-0.866,0.5){\circle*{0.15}}
             \bezier{90}(0.174,0.985)(-0.21,0.45)(-0.866,0.5) }
\def\slchel{\put(-0.14,0.985){\circle*{0.15}}
            \put(0.174,-0.985){\circle*{0.15}}
            \put(0.2,-0.985){\line(-1,6){0.325}} }
\def\slchtw{\put(0.174,-0.985){\circle*{0.15}}
            \put(-0.866,0.5){\circle*{0.15}}
            \bezier{90}(0.174,-0.985)(0,0)(-0.866,0.5) }
\def\slchtth{\put(0.866,0.5){\circle*{0.15}}
             \put(0.866,-0.5){\circle*{0.15}}
             \bezier{60}(0.866,0.5)(0.5,0)(0.866,-0.5) }
%
%
\def\achone{\put(0.174,0.985){\circle*{0.15}}
            \put(0.866,-0.5){\circle*{0.15}}
            \bezier{80}(0.174,0.985)(0.36,0.16)(0.866,-0.5) }
\def\achtwo{\put(0.174,-0.985){\circle*{0.15}}
            \put(0.866,0.5){\circle*{0.15}}
            \bezier{80}(0.174,-0.985)(0.36,-0.16)(0.866,0.5) }
\def\achthree{\put(-0.174,0.985){\circle*{0.15}}
            \put(-0.866,-0.5){\circle*{0.15}}
            \bezier{80}(-0.174,0.985)(-0.36,0.16)(-0.866,-0.5) }
\def\achfour{\put(-0.174,-0.985){\circle*{0.15}}
            \put(-0.866,0.5){\circle*{0.15}}
            \bezier{80}(-0.174,-0.985)(-0.36,-0.16)(-0.866,0.5) }
\def\achfive{\put(-0.174,-0.985){\circle*{0.15}}
              \put(-0.866,-0.5){\circle*{0.15}}
              \bezier{60}(-0.174,-0.985)(-0.36,-0.48)(-0.866,-0.5) }
\def\achsix{\put(0.174,-0.985){\circle*{0.15}}
            \put(-0.866,-0.5){\circle*{0.15}}
            \bezier{80}(0.174,-0.985)(-0.21,-0.45)(-0.866,-0.5) }
\def\achseven{\put(0.174,0.985){\circle*{0.15}}
             \put(0.174,-0.985){\circle*{0.15}}
             \bezier{100}(0.174,0.985)(-0.3,0)(0.174,-0.985) }
\def\acheight{\put(-0.866,0.5){\circle*{0.15}}
             \put(-0.866,-0.5){\circle*{0.15}}
             \bezier{60}(-0.866,0.5)(-0.5,0)(-0.866,-0.5) }
\def\achnine{\put(0.174,0.985){\circle*{0.15}}
             \put(-0.866,-0.5){\circle*{0.15}}
             \bezier{90}(0.174,0.985)(0,0)(-0.866,-0.5) }
%
%
\def\dslone{\bepi{\sctw{25}{25}{4}{4} \slchone \slchtwo \slchthree}}
\def\dsltwo{\bepi{\sctw{25}{25}{4}{4} \slchone \slchfour \slchthree}}
\def\dslthree{\bepi{\sctw{25}{25}{4}{4} \slchone \slchtwo \slchfive}}
\def\dslfour{\bepi{\sctw{25}{25}{4}{4} \slchone \slchfour \slchfive}}
\def\dslfive{\bepi{\sctw{25}{25}{4}{4} \slchsix \slchtth}}
\def\dslsix{\bepi{\sctw{25}{25}{4}{4} \slchseven \slcheight}}
\def\dslseven{\bepi{\sctw{4}{25}{4}{25} \slchone \slchnine \slchthree}}
\def\dsleight{\bepi{\sctw{4}{25}{4}{25} \slchone \slchten \slchthree}}
\def\dslnine{\bepi{\sctw{4}{25}{4}{25} \slchone \slchnine \slchfive}}
\def\dslten{\bepi{\sctw{4}{25}{4}{25} \slchone \slchten \slchfive}}
\def\dslel{\bepi{\sctw{4}{25}{4}{25} \slchel \put(-0.866,0.5){\circle*{0.15}}
    \put(0.866,-0.5){\circle*{0.15}} \put(0.866,-0.5){\line(-5,3){1.72}} }}
\def\dsltw{\bepi{\sctw{4}{25}{4}{25} \slchtw \slcheight}}
%
%
\def\aslone{\bepi{\sctw{25}{25}{4}{4} \slchone \slchseven \slcheight}}
\def\asltwo{\bepi{\sctw{25}{25}{4}{4} \slchone \slchseven \achone}}
\def\aslthree{\bepi{\sctw{25}{25}{4}{4} \slchone \achtwo \slcheight}}
\def\aslfour{\bepi{\sctw{25}{25}{4}{4} \slchone \achone \achtwo}}
\def\aslfive{\bepi{\sctw{25}{25}{4}{4} \slchtwo \slchthree}}
\def\asltse{\bepi{\sctw{25}{4}{25}{4} \slchone \achsix \slchtwo}}
\def\asltei{\bepi{\sctw{25}{4}{25}{4} \slchone \achsix \slchfour}}
\def\asltni{\bepi{\sctw{25}{4}{25}{4} \slchone \achfive \slchtwo}}
\def\asltte{\bepi{\sctw{25}{4}{25}{4} \slchone \achfive \slchfour}}
\def\asltto{\bepi{\sctw{25}{4}{25}{4} \achnine \slchseven}}
\def\aslttt{\bepi{\sctw{25}{4}{25}{4}\slchel\put(-0.866,-0.5){\circle*{0.15}}
   \put(0.866,0.5){\circle*{0.15}}\put(0.866,0.5){\line(-5,-3){1.72}} }}
%
%
\def\di#1{\raisebox{0pt}[20pt][23pt]{#1}}
\def\wdi#1{b^{\fsl(2)}\biggl({\raisebox{0pt}[20pt][30pt]{#1}}\biggr)}
\def\spwdi#1{W\biggl( \mbox{#1} \biggr)}
\def\slrel#1#2#3#4#5#6#7{
\begin{eqnarray*}
\makebox[2pt]{} \wdi{#1}-\wdi{#2}-\wdi{#3}+\wdi{#4} 
& \\
\hfill = \,\,\, \wdi{#5} - \wdi{#6} #7 \phantom{verylongstring} &
\end{eqnarray*}
}
\slrel{\dslone}{\dsltwo}{\dslthree}{\dslfour}{\dslfive}{\dslsix};
\slrel{\aslone}{\asltwo}{\aslthree}{\aslfour}{\dslfive}{\aslfive}.
\caption{The recurrence relation
	for the complete $\fsl(2)$ weight system}\label{7fig:31}
\end{figure}

The projection to the subspace of primitive elements along the subspace of
decomposable elements has the following form.

\begin{proposition}~{\rm \cite{L97}}
The projection~$\pi_n$ takes a chord diagram~$d$ with~$n$ chords to the linear combination
$$
\pi_n:d\mapsto d-1!\sum_{I\sqcup J=V(d)}d|_Id|_J+2!\sum_{I\sqcup J\sqcup K=V(d)}d|_Id|_Jd|_K-\dots,
$$
where the summations run over all unordered partitions of the set $V(d)$ of the chords in~$d$ into disjoint union
of nonempty subsets, and $d|_I$ is the chord subdiagram induced from~$d$ by the subset $I\subset V(d)$ of chords.
\end{proposition}

\begin{table}
\begin{center}
\begin{tabular}{|c|c|c|c|c|}
\hline
$k$ & $d\in\cM_{2k}$& $b^{\fsl_2}(d)$  & \phantom{di}$\gamma(d)$\phantom{di} & $R_k(d)$\\
&&$b^{\fsl_2}(\pi_{2k}(d)) $&&\\
\hline
$2$&\begin{picture}(20,30)(0,-2)
\put(0,10){\circle{24}}
\put(0,-2){\line(0,1){24}}
\put(-8,2){\line(1,1){16}}
\put(8,2){\line(-1,1){16}}
\put(-12,12){\line(2,-3){9}}
\end{picture}&$c(c-1)^2(c-2)$&\begin{picture}(20,30)(0,0)
\put(0,10){\circle*{3}}
\put(10,0){\circle*{3}}
\put(10,20){\circle*{3}}
\put(20,20){\circle*{3}}
\put(0,10){\line(1,-1){10}}
\put(0,10){\line(1,1){10}}
\put(10,20){\line(0,-1){20}}
\put(10,20){\line(1,0){10}}
\end{picture}&$0$\\
&&$-2c$&&\\
\hline
$2$&\begin{picture}(20,30)(0,-2)
\put(0,10){\circle{24}}
\put(-10,4){\line(1,0){20}}
\put(-8,2){\line(0,1){16}}
\put(-10,16){\line(1,0){20}}
\put(8,2){\line(0,1){16}}
\end{picture}&$c^4-4c^3+8c^2-4c$&\begin{picture}(20,30)(0,0)
\put(0,0){\circle*{3}}
\put(20,0){\circle*{3}}
\put(0,20){\circle*{3}}
\put(20,20){\circle*{3}}
\put(0,0){\line(1,0){20}}
\put(0,0){\line(0,1){20}}
\put(20,20){\line(0,-1){20}}
\put(20,20){\line(-1,0){20}}
\end{picture}&$1$\\
&&$2c^2-4c$&&\\
\hline
$2$&\begin{picture}(20,30)(0,-2)
\put(0,10){\circle{24}}
\put(-10,4){\line(1,0){20}}
\put(-8,2){\line(1,1){16}}
\put(-10,16){\line(1,0){20}}
\put(8,2){\line(-1,1){16}}
\end{picture}&$c^4-5c^3+10c^2-5c$&\begin{picture}(20,30)(0,0)
\put(0,0){\circle*{3}}
\put(20,0){\circle*{3}}
\put(0,20){\circle*{3}}
\put(20,20){\circle*{3}}
\put(0,0){\line(1,0){20}}
\put(0,0){\line(0,1){20}}
\put(0,0){\line(1,1){20}}
\put(20,20){\line(0,-1){20}}
\put(20,20){\line(-1,0){20}}
\end{picture}&$1$\\
&&$2c^2-5c$&&\\
\hline
$2$&\begin{picture}(20,30)(0,-2)
\put(0,10){\circle{24}}
\put(0,-2){\line(0,1){24}}
\put(-12,10){\line(1,0){24}}
\put(-8,2){\line(1,1){16}}
\put(8,2){\line(-1,1){16}}
\end{picture}&$c^4-6c^3+13c^2-7c$&\begin{picture}(20,30)(0,0)
\put(0,0){\circle*{3}}
\put(20,0){\circle*{3}}
\put(0,20){\circle*{3}}
\put(20,20){\circle*{3}}
\put(0,0){\line(1,0){20}}
\put(0,0){\line(0,1){20}}
\put(0,0){\line(1,1){20}}
\put(20,20){\line(0,-1){20}}
\put(20,20){\line(-1,0){20}}
\put(0,20){\line(1,-1){20}}
\end{picture}&$1$\\
&&$2c^2-7c$&&\\
\hline
$3$&\begin{picture}(20,30)(0,-2)
\put(0,10){\circle{24}}
\put(-10,16){\line(1,0){20}}
\put(-10,4){\line(1,0){20}}
\put(-12,8){\line(2,3){9}}
\put(-12,12){\line(2,-3){9}}
\put(12,8){\line(-2,3){9}}
\put(12,12){\line(-2,-3){9}}
\end{picture}&$c^6-6c^5+15c^4-18c^3+18c^2-8c$&\begin{picture}(20,30)(0,0)
\put(3,0){\circle*{3}}
\put(21,0){\circle*{3}}
\put(3,24){\circle*{3}}
\put(21,24){\circle*{3}}
\put(-3,12){\circle*{3}}
\put(27,12){\circle*{3}}
\put(3,0){\line(1,0){18}}
\put(3,24){\line(1,0){18}}
\put(-3,12){\line(1,2){6}}
\put(-3,12){\line(1,-2){6}}
\put(27,12){\line(-1,2){6}}
\put(27,12){\line(-1,-2){6}}
\end{picture}&$1$\\
&&$2c^3+3c^2-8c$&&\\
\hline
$3$&\begin{picture}(20,30)(0,-2)
\put(0,10){\circle{24}}
\put(-10,16){\line(1,0){20}}
\put(-10,4){\line(1,0){20}}
\put(-12,8){\line(2,3){9}}
\put(-12,12){\line(2,-3){9}}
\put(2,-2){\line(1,3){7}}
\put(2,22){\line(1,-3){7}}
\end{picture}&$c^6-8c^5+28c^4-57c^3+66c^2-27c$&\begin{picture}(20,30)(0,0)
\put(3,0){\circle*{3}}
\put(21,0){\circle*{3}}
\put(3,24){\circle*{3}}
\put(21,24){\circle*{3}}
\put(-3,12){\circle*{3}}
\put(27,12){\circle*{3}}
\put(3,0){\line(1,0){18}}
\put(3,24){\line(1,0){18}}
\put(-3,12){\line(1,2){6}}
\put(-3,12){\line(1,-2){6}}
\put(27,12){\line(-1,2){6}}
\put(27,12){\line(-1,-2){6}}
\put(-3,13){\line(5,2){24}}
\put(27,13){\line(-5,2){24}}
\end{picture}&$0$\\
&&$20c^2-27c$&&\\
\hline
$3$&\begin{picture}(20,30)(0,-2)
\put(0,10){\circle{24}}
\put(0,-2){\line(0,1){24}}
\put(-12,10){\line(1,0){24}}
\put(-10,16){\line(2,-1){21}}
\put(-6,21){\line(1,-2){11}}
\put(10,16){\line(-2,-1){21}}
\put(6,21){\line(-1,-2){11}}
\end{picture}&$c^6-15c^5+115c^4-430c^3+657c^2-295c$&\begin{picture}(20,30)(0,0)
\put(3,0){\circle*{3}}
\put(21,0){\circle*{3}}
\put(3,24){\circle*{3}}
\put(21,24){\circle*{3}}
\put(-3,12){\circle*{3}}
\put(27,12){\circle*{3}}
\put(3,0){\line(1,0){18}}
\put(3,24){\line(1,0){18}}
\put(-3,12){\line(1,2){6}}
\put(-3,12){\line(1,-2){6}}
\put(27,12){\line(-1,2){6}}
\put(27,12){\line(-1,-2){6}}
\put(-3,13){\line(5,2){24}}
\put(27,13){\line(-5,2){24}}
\put(3,0){\line(0,1){24}}
\put(3,1){\line(4,5){19}}
\put(3,1){\line(5,2){24}}
\put(21,0){\line(0,1){24}}
\put(21,1){\line(-4,5){19}}
\put(21,1){\line(-5,2){24}}
\put(-3,12){\line(1,0){30}}
\end{picture}&$8$\\
&&$16c^3-284c^2+295c$&&\\
\hline
\end{tabular}
\end{center}
\caption{Values of the weight systems $b^{\fsl_2}(\cdot)$, $b^{\fsl_2}(\pi(\cdot))$, $R_k$ on certain chord diagrams}
\end{table}

Here is more evidence supporting the conjecture:
\begin{enumerate}
\item the weight system $b^{\fsl_2}$,
similarly to the weight system~$R_k$, depends on the intersection graphs of the chord diagram
rather than on the diagram itself---see~\cite{CL07};
\item the value~$R_k(d_1d_2)$ on a product of two nontrivial diagrams with
$2k$ chords in total is~$0$, since there are no $2k$-gons in the intersection graph;
on the other hand, $\pi_{2k}(d_1d_2)=0$, since the diagram $d_1d_2$ is decomposable;
\item for a chord diagram~$d$ with~$2k$ chords having a leaf (a chord intersecting
only one other chord), $R_k(d)=0$; it can also be easily proved that
the degree of the polynomial $b^{sl_2}(\pi_{2k}(d))$ is less than~$k$.
\end{enumerate}

\begin{remark}
The same argument as in the proof of Proposition~\ref{pge} allows one to extend the value of the
$b^{\fsl_2}$ weight system to a $4$-invariant of graphs with up to~$6$ edges. Namely,
we can set the value of this invariant on the $5$-wheel to be equal to $c^6-10c^5+50c^4-139c^3+176c^2-72c$,
with the projection to the subspace of primitive elements $-6c^3+70c^2-72c$,
and on the $3$-prism to be equal to $c^6-9c^5+40c^4-108c^3+146c^2-63c$,
with the projection to the subspace of primitive elements $-2c^3+58c^2-63c$.
\end{remark}

\section{Note added in proof}\label{s5}

After the first version of the present text has been spread as a preprint and submitted for publication,
a paper by Bar-Natan and Vo~\cite{BNV14} appeared. In this paper,
Conjecture~\ref{conj-sl2-hc} is shown to be true, and now it can be stated as a theorem:

\begin{theorem}{\rm\bf \cite{BNV14}}\label{th-sl2-hc}
For any chord diagram~$d$ with~$2k$ chords, the coefficient of~$c^k$ in the value
$b^{\fsl_2}(\pi_{2k}(d))$ is $2R_k(d)$.
\end{theorem}

In fact, Bar-Natan and Vo show that the proof essentially constitutes a part of
the proof of the Melvin--Morton--Rozansky Conjecture in~\cite{BNG96}.

One of the key ingredients in the proof is Proposition 3.13 in~\cite{BNG96}
which reads as follows. Denote by $w_C$ the $\Z$-valued function on chord diagrams
equal to~$1$ if the adjacency matrix of the intersection graph of the diagram
is nondegenerate over~$\Z/2\Z$, and equal to~$0$ otherwise. The function~$w_C$ is extended to~$\cM$
by linearity. It is easy to show that~$w_C$
satisfies the so-called $2$-term relation, $w_C(d)=w_C(\tilde d_{AB})$, for
any chord diagram~$d$ and any pair of chords~$A$ and~$B$ with neighboring ends in it~\cite{BNG96};
the $4$-term relation is an obvious corollary of the $2$-term one.

\begin{proposition}~{\rm\bf\cite{BNG96}}
$$
R_k(d)=-w_C(\pi_{2k}(d))
$$
for any chord diagram~$d$ with $2k$~chords.
\end{proposition}

This argument allows one to prove the following statement generalizing Proposition~\ref{pge}
to arbitrary values of~$k$.

\begin{theorem}
For arbitrary~$k$, the function~$R_k$ extends to a $4$-invariant of graphs.
\end{theorem}

Indeed, the function~$w_C$
is known to be extendable to a multiplicative $\Z$-valued $4$-invariant of graphs, see~\cite{CL07} or~\cite{M03}:
set~$w_C(\Gamma)$ to be equal to~$1$ if the adjacency matrix of~$\G$
is nondegenerate over~$\Z/2\Z$, and equal to~$0$ otherwise.
It satisfies obviously the $2$-term relation for graphs, $w_C(\G)=w_C(\widetilde\G_{AB})$,
for any graph~$\G$ and any pair of vertices~$A,B$ in it.
For intersection graphs with $2k$~vertices,
the projection of this graph invariant to primitive elements coincides with the invariant~$-R_k$.
Hence, the projection to primitive elements $w_C(\pi_{2k}(\G))$ of this $4$-invariant on arbitrary graph~$\G$
with~$2k$ vertices is a $4$-invariant that coincides with~$-R_k(\G)$ if~$\G$
is an intersection graph. \QED


\vspace{1cm}
\noindent
National Research University Higher School of Economics\\
\noindent
7 Vavilova Moscow 117312 Russia\\

\vspace{1cm}

\noindent
\begin{tabular}{ll}
E.~Kulakova &lenaetopena@gmail.com\\
S.~Lando &lando@hse.ru\\
T.~Mukhutdinova &kassalanche@gmail.com\\
G.~Rybnikov &grigory.rybnikov@gmail.com
\end{tabular}

\begin{thebibliography}{99}

\bibitem{BN95}   {D.~Bar-Natan,}
    	{\it On Vassiliev knot invariants},
        Topology, vol.~{\bf 34}, no.~2 (1995), 423--472.


\bibitem{BNG96} D.~Bar-Natan, S.~Garoufalidis, {\it On the Melvin--Morton--Rozansky Conjecture}, Inventiones Mathematicae, vol.~{\bf 125} (1996), 103--133.

\bibitem{BNV14} D.~Bar-Natan, H.~Vo, {\it Proof of a conjecture of Kulakova et al. related to the $\fsl_2$ weight system}, preprint, arXiv:1401.0754 [math.QA] (2014).

\bibitem{CDM12} S.~Chmutov, S.~Duzhin, Y.~Mostovoy, {\it Introduction to Vassiliev Knot Invariants},
Cambridge University Press, 2012, ISBN 978-1-107-02083-2.
\verb!http://www.pdmi.ras.ru/~duzhin!

\bibitem{CL07} S.~V.~Chmutov, S.~K.~Lando,
{\it Mutant knots and intersection graphs},
Algebraic and Geometric Topology, vol.~{\bf 7} (2007), 101--120.

\bibitem{CV97} 	{S.~V.~Chmutov, A.~N.~Varchenko,}
        {\it Remarks on the Vassiliev knot invariants
        coming from $\fsl_2$},
        Topology, vol.~{\bf 36} (1997), 153--178.

\bibitem{K93}	{M.~Kontsevich,}
        {\it Vassiliev knot invariants},
        in: Adv. in Soviet Math., vol.~{\bf 16} (1993), part 2, 137--150.


\bibitem{L97} 	{S.~K.~Lando},
    	{\it On primitive elements in the bialgebra of
    	chord diagrams},
    	in: Amer. Math. Soc. Transl. Ser.~$2$,
    	AMS, Providence RI, 1997, vol.~{\bf 180}, 167--174.

\bibitem{L00}      {S.~K.~Lando,}
    	{\it On a Hopf algebra in graph theory},
    	J. Comb. Theory, Ser. B, vol.~{\bf 80} (2000), 104--121.

\bibitem{LZ04} S.~Lando, A.~Zvonkin, {\it Graphs on surfaces and their applications}
(Chapter~6), Springer, 2004.

\bibitem{M03}    B.~Mellor, {\it A few weight systems arising from intersection graphs},
    Michigan Math. J., vol.~{\bf 51}, no. 3 (2003), 509--536.

\end{thebibliography}
\end{document}